\documentclass{amsart}

 \newcommand\Pn[1]{{\bf P}^{#1}} 
 \newcommand\Pnd[1]{{\check{\bf P}}^{#1}}
 \newtheorem{theorem}{{\rm T\sc heorem}}[section] 
 \newtheorem{lemma}[theorem]{{\rm L\sc emma}} 
  
 \newtheorem{proposition}[theorem]{{\rm P\sc roposition}}

 \newtheorem{problem}[theorem]{Problem} 
  
 \newtheorem{re}[theorem]{{\rm R\sc emark}}

 \newtheorem{alemma}[theorem]{{\rm A\sc polarity lemma}}
  
 \newcommand\C{{\bf C}} 
 \newcommand\Gr{{\bf G}}
 
 \newcommand\Hom{\mathop{\rm Hom}\nolimits}

 \newcommand\rto{--->}

\begin{document}
 \title {Canonical curves and varieties of sums of powers of cubic 
 polynomials}
 \author{Atanas Iliev and Kristian Ranestad}
 
 \begin{abstract} 
In this note we show that the apolar cubic forms associated to 
 codimension 2 linear sections of canonical curves of genus  
 $g\geq 11$ 
 are special with respect to their presentation as sums of cubes. 
 \end{abstract}
 \maketitle 
 {\footnotetext{Mathematics Subject Classification 13H10 (Primary), 
 13E10,
 14M05, 
 14H45 (Secondary)}}
 \pagestyle{myheadings}\markboth{\textsc{Atanas Iliev and Kristian 
 Ranestad}}{\textsc{Canonical curves and varieties of sums of powers of cubic 
 polynomials}}
 \section{Introduction}{\label{s1}}

 In a graded Artinian Gorenstein ring $A$ with socle degree $d$, 
 multiplication defines (up to scalar) a homogeneous form $f$ of 
 degree $d$, called the socle degree generator, dual polynomial or 
 apolar polynomial of $A$.  Codimension $2$ linear 
sections of a canonical curve of genus $g$ define Artinian Gorenstein quotients 
 of the homogeneous coordinate ring of the 
curve.  These quotients have socledegree $3$ and therefore define (up to 
scalar) cubic forms in $g-2$ variables.
A dimension count shows that a general 
cubic form is not obtained this way when $g\geq 8$. 
 While a general cubic form in $g-2$ variables cannot be written as a 
 sum of less than ${\frac 16}g(g-1)$ cubes, our main result says that the 
cubic forms apolar to a general codimension $2$ linear 
section of a general canonical curve of genus $g\geq 11$  can be 
written as a sum of $2g-4$ cubes. 

Our methods give results concerning the variety of different 
powersum presentations.  In particular we obtain partial results for 
genus $g=9$ (cf. \ref{s3}).  Results for $g\leq 6$ are classical, while $g=7$ and 
$g=8$ was treated in \cite{RS} and \cite{IR}. 

Powersum presentations of forms from a more algebraic viewpoint have 
been studied extensively in \cite{IK}.  

We work throughout over the complex numbers $\C$.

\subsection{Powersum presentations}
  Let $f \in \C[x_0,\ldots,x_n]$ be a homogeneous form of degree d, 
  then $f$ can be written as a sum
of powers of linear forms
$$f = l_1^d+ \ldots + l_s^d$$
for $s$ sufficiently large. Indeed, if we identify the map $l \mapsto l^d$ with the $d^{th}$ Veronese
embedding $\Pn n \hookrightarrow \Pn {N_d}$, where $N_d= 
{{n+d}\choose n}-1$, this amounts to say 
that the image spans $\Pn {N_d}$. Fixing $(d,n)$, the minimal number $s$ of summands 
needed varies with $f$, of course.  A simple dimension count shows that
$$ s \ge \lceil {\frac 1{n+1}}{n+d \choose n} \rceil $$
for a {\bf general} $f$. With a few exceptions equality holds by a result of Alexander
and Hirschowitz \cite{AH} combined with Terracini's Lemma (cf. 
\cite{Iar}) :
\begin {theorem}{\label{1.1}}{\rm (Alexander,Hirschowitz)}  A general form $f$ of degree $d$ in $n+1$ variables
is a sum of $\lceil {\frac 1{n+1}}{n+d \choose n} \rceil$ powers of linear 
forms, unless\vskip 1pt 
 $d=2$, where $s=n+1$ instead of $\lceil {\frac {n+2}2}\rceil$, or\vskip 1pt 
 $d=4$ and $n=2,3,4$, where $s=6,10,15$ instead of $5,9,14$ respectively, or\vskip 1pt 
$d=3$ and $n=4$, where $s=8$ instead of $7$.\end{theorem}

Let $F=Z(f)\subset \Pn n$ be the hypersurface defined by $f$.  For a 
linear form $l$ we denote by $L$ 
 the point in $\Pnd n$ of the 
 hyperplane $Z(l)\subset\Pn n$.  Then we 
define, as in \cite{RS}, 
the {\bf v}ariety 
 of {\bf s}ums of {\bf p}owers 
 as the closure 
 $$VSP(F,s) = 
 \overline{\{ \{L_{1},\ldots,L_{s} \} \in Hilb_s(\Pnd n) \mid \exists 
 \lambda_i \in \C 
 : f=\lambda_1 l_1^d+\ldots+\lambda_s l_s^d \} }$$ 
 of the set of powersums presenting $f$ in the Hilbert scheme (cf. 
 \cite{RS}).  Notice that taking $d^{\rm 
th}$ roots of the $\lambda_{i}$, we can put them into the forms $l_{i}$.   
We study these varieties of sums of powers 
 using apolarity.

 \subsection{Apolarity}{\label{apolarity}}(cf. \cite{RS}). 
 Consider $R=\C[x_0,\ldots,x_n]$ and 
 $T=\C[\partial_0,\ldots,\partial_n]$. $T$ acts on $R$ by 
 differentiation: 
 $$\partial^{\alpha}\cdot x^{\beta} = \alpha!{\binom{\beta}{\alpha} } 
 x^{\beta-\alpha}$$ 
 if $\beta \geq \alpha$ and 0 otherwise. Here $\alpha$ and $\beta$ are 
 multi-indices, 
 ${\binom{\beta}{\alpha}} = \prod {\binom{\beta_i}{\alpha_i}}$ and so on.
 One can interchange the role of R and T by defining 
 $$x^{\beta}\cdot \partial^{\alpha} = \beta!{\binom{\beta}{\alpha} } 
 \partial^{\alpha-\beta}.$$ 
 This action defines a perfect pairing between forms of degree d and 
 homogeneous differential operators of order d. In particular, 
 $R_{1}$ and $T_{1}$ are natural dual vector spaces.  Therefore the projective spaces with 
 coordinate ring $R$ and $T$ respectively are natural dual to each 
 other, we denote them by $\Pn n$ and $\Pnd n$.  
 A point $a=(a_0,\ldots,a_n) \in {\Pnd n}$ defines a form $l_a = \sum 
 a_i x_i \in R_{1}$, and for 
 a form $D\in T_{e}$ 
 $$D\cdot l_a^d=e!{\binom{d}{e} }D(a)l_a^{d-e},$$ 
 when $e\leq d$.
 In particular $$D\cdot l_a^d=0 \iff D(a)=0 \leqno(*)$$ if $e \leq d$. 
 More generally we say that homogeneous forms $f\in R$ and $D\in T$ 
 are {\bf apolar} if $f\cdot D=D\cdot f=0$ (According to 
 Salmon (1885) \cite{Sal} the term was coined by Reye).
 
 Apolarity allows us to associate  an Artinian Gorenstein 
 graded quotient ring of $T$ to a form: For $f\in R$ a homogeneous form of degree $d$ and $F =Z(f) \subset 
 \Pn n$ define 
 $$ F^{\bot} = f^{\bot} = \{D \in T | D\cdot f=0 \}$$ 
 and $$A^F = T/F^{\bot}.$$
 The socledegree of $A^F$ is $d$, since 
 $$D^{\prime}\cdot (D\cdot f) = 0 \hskip3pt \forall D^{\prime} \in T_1 \iff D\cdot f = 0 \hskip3pt 
 or \hskip3pt D \in T_d.$$ 
 In particular the socle of $A^F$ is 1-dimensional, and $A^F$ is 
  Gorenstein.  It is called the  apolar Artinian 
 Gorenstein ring of $F$ .
 Conversely for a graded Gorenstein ring $A = T/I$ with socledegree $d$, 
 multiplication in A induces a linear form $f\colon {\rm Sym}_d(T_1) \to \C$ 
 which can be identified with a homogeneous 
 polynomial $f \in R$ of degree $d$. This proves: 
 \setcounter{theorem}{1} 
 \begin{lemma}{\label{1.2}} {\rm (Macaulay, \cite{Mac})} 
 The map $F \mapsto A^F$ is a bijection between hypersurfaces $F=Z(f) 
 \subset \Pn n$ of 
 degree d and graded Artinian Gorenstein quotient rings $A = T/I$ of T 
 with socledegree d.\end{lemma}
 
 Let $X\subset \Pn {n+m+1}$ be a $m$-dimensional arithmetic Gorenstein 
 variety.  Let $S(X)$ be the homogeneous coordinate ring of $X$, and 
 let $h_{1},\ldots,h_{m+1}$ be general linear forms and set $L=Z(h_{1},\ldots,h_{m+1})$. Then by definition 
 $S(X)/(h_{1},\ldots,h_{m+1})$ is Artinian Gorenstein, i.e. by Macaulays 
 result the apolar Artinian Gorenstein ring of a $(n-1)$-dimensional hypersurface $F_{L}$ of degree $d$, 
 the socledegree of the ring.   $L$ is a linear space of 
 dimension $n$ and by apolarity $F_{L}=Z(f_{L})$ is a hypersurface in the
 dual space to $L$. We say that $F_{L}$ is apolar to the (empty) linear section $L\cap X$. Hence, there is a rational map
 $$\alpha_{X}:\Gr (n+1,m+n+2)\rto H_{n,d}$$
 Where $H_{n,d}$ is the space of $(n-1)$-dimensional hypersurfaces of degree $d$  
 modulo the action of 
$PGL(n+1,k)$.

A canonical curve  $C\subset \Pn {g(C)-1}$ is arithmetic Gorenstein, 
i.e. the homogeneous coordinate ring $S(C)$ is Gorenstein. Let $h_{1}, 
h_{2}\in S(C)$ be two general linear 
forms, then the quotient 
$S(C)/(h_{1},h_{2})$ is Artinian Gorenstein with values of the Hilbert 
function:  $1,g-2,g-2,1$.  Its socledegree 
is therefore $3$.  Thus we obtain a map
$$\alpha_{C}:\Gr (g(C)-2,g(C))\rto H_{g(C)-3,3}$$
 to the space of cubic hypersurfaces of dimension $g(C)-4$.  
 We shall study the image of this map. In particular 
 we shall study the variety of sums of powers of the cubic 
 hypersurfaces in this image.

 \subsection {Variety of apolar subschemes} 

  Let $F=Z(f) \subset \Pn n$ 
 denote a hypersurface of degree d. We call a subscheme 
 $\Gamma \subset \Pnd n$ 
 $\bf {apolar} $ to $F$, if the homogeneous ideal $I_{\Gamma} \subset 
 F^{\bot} \subset T$.  
 
 \begin{alemma}{\label{1.4}} Let $l_1,\ldots, l_s$ be linear forms in $R$, 
 and let $L_{i} \in \Pnd n$ be 
 the corresponding points in the dual space. Then 
 $f = \lambda_{1}l_1^d+\ldots+\lambda_{s}l_s^d$ for some 
 $\lambda_{i}\in \C^{*}$ if and only if $\Gamma = \{L_{1},\ldots 
 ,L_{s}\} \subset \Pnd n$ is apolar to $F=Z(f)$.\end{alemma}
 \proof 
 Assume $f = \lambda_{1}l_1^d+\ldots+\lambda_{s}l_s^d$. If $g\in 
 I_{\Gamma}$, then $g\cdot l_i^d$=0 for all $i$ by (*), so by 
 linearity $g\in F^{\bot}$. Therefore $\Gamma$ is apolar to $F$. \par 
 For the converse, assume that $I_{\Gamma}\subset F^{\bot}$. Then we have 
 surjective maps between 
 the corresponding homogeneous coordinate rings 
 $$T\to A_{\Gamma}=T/I_{\Gamma}\to A^F.$$ 
 Consider the dual inclusions of the degree $d$ part of these rings:  
 $$\Hom (A^F_d,\C)\to \Hom ((A_{\Gamma})_d,\C)\to \Hom (T_d,\C).$$ 
 $D\mapsto D\cdot f$ generates the first of these spaces, while the second is 
 spanned by the forms 
 $D\mapsto D\cdot l_i^d$. Thus $f$ lies in the span of the $l_i^d$. \qed
 
 \vskip 2pt
This is the crucial lemma in the study of powersum presentations of 
$f$.  Furthermore it allows us to define a {\bf v}ariety of a{\bf p}olar 
{\bf s}ubschemes 
to $f$, which naturally extends our definition of the variety of sums 
of powers. 
 $$VPS(F,s) = 
 \overline{\{\Gamma\in Hilb_s(\Pnd n) \mid I_{\Gamma}\subset 
 F^{\bot}\} },$$ 
where $Hilb_s(\Pnd n)$ is the Hilbert scheme of length $s$ subschemes 
 of $\Pnd n$.  Clearly $VSP(F,s)$ is the closure of the set 
 parametrizing smooth subscemes in $VPS(F,s)$. In 
 general they do not coinside.

  \section {Apolar varieties of singular sections}{\label{s2}}
  \subsection{Apolar varieties}{\label{apolar}}
  Let $X\subset \Pn {n+m+1}$ be a reduced and irreducible 
  $m$-dimen-\break 
  sional nondegenerate variety of degree $d\geq 3$ and codimension 
  $n+1\geq 2$. Let $p\in X$ be a general smooth point.
  Let $C_{p}X$ be the cone over $X$ with vertex at $p$.  Since $p$ is 
  a smooth point, the degree of the cone $C_{p}X$ is $d-1$, while the 
  dimension is $m+1$.  Clearly $X\subset C_{p}X$.  
 
  We apply this simple construction to describe powersum presentations of 
hypersurfaces in the image of the map $\alpha_{X}$ in \ref{1.2}. 
Let again $X\subset \Pn {n+m+1}$ be a $m$-dimensional arithmetic Gorenstein 
 variety of degree $d$.  
 Fix a general $n$-dimensional linear subspace 
 $L\subset \Pn {n+m+1}$, in particular we fix the hypersurface $F_{L}$ in the image of $\alpha_{X}$.
  Let $p$ be a smooth point on $X$, then the intersection $C_{p}X\cap L$ is clearly nonempty, and if it is proper it 
 is $0$-dimensional of degree $d-1$.  We may assume that this 
 intersection is proper and smooth for a general $L$ or general $p$, so 
 we get an apolar subscheme of degree $d-1$ to $F_{L}$, i.e. a point in 
 $VSP(F_{L},d-1)$.  We have shown:
  \begin{proposition}{\label{2.1}} Let $X\subset \Pn {n+m+1}$ be a $m$-dimensional arithmetic Gorenstein 
 variety of degree $d$, and let $L\subset \Pn {n+m+1}$ be a $n$-dimensional 
 linear subspace such that $L\cap X=\emptyset$.  Let $F_{L}$ be the 
 associated apolar hypersurface. Then there is a rational map  $X\rto 
 VSP(F_{L},d-1)$ defined by $p\mapsto C_{p}X\cap L$.
\end{proposition}
 
 \begin{problem} When is this map a morphism?  When can $F_{L}$ and 
 $X$ be recovered from the image of this map?\end{problem}
 
 We may improve slightly on the degree of the apolar subschemes by 
 considering cones on special linear sections of $X$.\par
\subsection{Tangent hyperplane sections}{\label{tangent}}
Let $X\subset \Pn {n+m+1}$ be a reduced and irreducible 
  $m$-dimensional nondegenerate variety of degree $d$ and codimension 
  $n+1\geq 2$.   We assume additionally that $X$ satisfies the 
  following condition: 
 \begin{tabbing} \hspace{1.5cm} \= 
  A general tangent hyperplane 
  section of $X$ has a double point at  \\
  $({\ast\ast})$ \>
  the point of tangency, and the projection of the tangent hyperplane \\ 
  \> section from the point of tangency is 
  birational.\\
  \end{tabbing}

  In particular, $X$ is not a scroll and $d\geq 4$. 
  Let $p\in X$ be a general smooth point.
  Let $H_{p}$ be a general hyperplane tangent to $X$ at $p$. Since 
 $H_{p}\cap X$ has multiplisity $2$ at $p$ and the 
 projection of $H_{p}\cap X$ from $p$ is birational, the image of the 
 projection is $(m-1)$-dimensional of degree $d-2$.   Therefore 
 $H_{p}\cap X$ is contained in an $m$-dimensional cone $C_{p}(H_{p}\cap 
 X)$ of degree $d-2$ 
 with vertex at $p$.  Similarly, if $H_{p}$ and $H_{p}'$ are two 
 general hyperplanes tangent at $p$, then the intersection $H_{p}\cap H_{p}'\cap X$
 has a singularity at $p$ of multiplisity 4, the complete 
 intersection of two singularities of multiplisity $2$. 
 In this case we say that the codimension $2$ space $H_{p}\cap H_{p}'$ 
 is {\bf doubly tangent} to $X$ at $p$.  If
 
 \begin{tabbing} $({\ast\ast\ast})$\hspace{1.2cm} \= 
  the projection of $H_{p}\cap H_{p}'\cap X$ 
 from $p$  is birational,  \\
   \end{tabbing}
   then the image $(m-2)$-dimensional of degree $4$ less 
 than the degree of $X$.  Hence $H_{p}\cap H_{p}'\cap X$ is contained in 
 a $(m-1)$-dimensional cone $C_{p}(H_{p}\cap H_{p}'\cap X)$ of degree $d-4$ with vertex at $p$. 
 This proves the
 
 \begin{lemma}{\label{2.2}} Let $X\subset \Pn {n+m+1}$ be a smooth 
 $m$-dimensional nondegenerate variety of degree $d$, and assume 
 that $X$ satisfies condition $(\ast\ast)$. Let $p\in X$ be a general smooth 
 point, and let $H_{p}$ be a general 
 hyperplane tangent to $X$ at $p$.  Then the cone $C_{p}(H_{p}\cap 
 X)$ is an $m$-dimensional variety of degree $d-2$ that contains $H_{p}\cap 
 X$.
 Assume furthermore 
 that $X$ satisfies condition $(\ast\ast\ast)$, and let $H_{p}$ and $H_{p}'$ be two general 
 hyperplanes tangent to $X$ at $p$.  Then the cone $C_{p}(H_{p}\cap H_{p}'\cap X)$ 
 is a $(m-1)$-dimensional variety of degree $d-4$ that contains $H_{p}\cap H_{p}'\cap X$.
 \end{lemma}
 
As above we apply this lemma to describe powersum presentations of 
hypersurfaces in the image of the map $\alpha_{X}$ in section 
\ref{apolarity}. 
Let again $X\subset \Pn {n+m+1}$ be a $m$-dimensional arithmetic Gorenstein 
 variety of degree $d$, with $m\geq 1$.  We assume additionally 
 that $X$ satisfies condition $(\ast\ast)$.
 Fix a general $n$-dimensional linear subspace 
 $L\subset \Pn {n+m+1}$, in particular we fix the hypersurface $F_{L}$ in the image of $\alpha_{X}$.
  If $L\subset H_{p}$ where $H_{p}$ is a 
 general hyperplane tangent at $p$, then according to \ref{2.1} 
 there is a $(m-1)$-dimensional variety $Y\supset H_{p}\cap X$ of degree 
 $d-2$.
 The intersection $Y\cap L$ is clearly nonempty, and if it is proper it 
 is $0$-dimensional of degree $d-2$.  We may assume that this 
 intersection is proper for a general $L$, so we get a point in 
 $VSP(F_{L},d-2)$.
 Let $\check {X}\subset \Pnd {m+n+1}$ be the dual variety of $X$, i.e. 
 the set of hyperplanes tangent 
 at some point $p\in X$.
 Then we have set up a rational map
 $$\check {X}_{L}\to VSP(F_{L},d-4)$$
 where $\check {X}_{L}=\{[H]\in \check {X}_{L}|H\supset L\}$.
 The subvariety 
 $\check {X}_{L}$ has dimension $m-{\rm codim}\check {X}$, which 
 equals $m-1$ when $\check {X}$ is nondegenerate.  In particular this is the case 
 when $X$ is a curve.
\par
 Similarly, assume that $X$ satisfies $(\ast\ast\ast)$, and  let $L\subset H_{p}\cap H_{p}'$ where $H_{p}$ and $H_{p}'$ are two 
 general hyperplanes tangent at $p$. Then, according to \ref{2.1}, 
 there is a $(m-1)$-dimensional variety 
 $$Y\supset H_{p}\cap H_{p}'\cap X$$
 of degree $d-4$.
 The intersection $Y\cap L$ is clearly nonempty, and if it is proper it 
 is $0$-dimensional of degree $d-4$.  We may assume that this 
 intersection is proper for a general $L$, so we get a point in $VSP(F_{L},d-4)$.
 Let $Z_{X}\subset \Gr (m+n,m+n+2)$ be the set of codimension 
 $2$ subspaces doubly tangent at some point $p\in X$.
 Then we have set up a rational map
 $$\Gr (m+n, m+n+2)\supset Z_{L}\rto VSP(F_{L},d-4)$$
 where $Z_{L}=\{[V]\in Z_{X}|V\supset L\}$.
 If the dual variety of $X$ is nondegenerate, then the subvariety 
 $Z_{X}$ has dimension $m+2(n-1)$.  The codimension in 
 $\Gr (m+n, m+n+2)$ of subspaces that contains $L$ is 
 $2(m+n)-2(m-1)=2n+2$  so the expected dimension of $Z_{L}$ is $m-4$.
 \par
 Notice that it is essential for the dimension count that $X$ is not a cone, i.e. that the 
 dual variety is nondegenerate.
 
 \begin{proposition}{\label{2.4}} Let $X\subset \Pn {n+m+1}$ be a $m$-dimensional arithmetic Gorenstein 
 variety of degree $d$, with $m\geq 1$.  Assume that $X$ satisfies 
 the condition $(\ast\ast)$ and has nondegenerate dual variety.    Let 
 $L\subset \Pn {n+m+1}$ be a general $n$-dimensional linear subspace, 
 and let $F_{L}=\alpha_{X}([L])$ be the hypersurface apolar to $L\cap 
 X$. Then $VSP(F_{L},d-2)\not=\emptyset$ and of dimension at least $m-1$.
 Assume furthermore that $m\geq 4$ and that $X$ also satisfies condition 
 $(\ast\ast\ast)$.   Then  
 the dimension of $VSP(F_{L},d-4)$ is at least $m-4$, when $m\geq 4$ 
 and $L$ is contained in at least one codimension $2$ linear space doubly tangent 
 to $X$.
 \end{proposition}
 
 \section {Canonical curves and apolar cubic polynomials}{\label{s3}}
 For a general cubic $n$-fold $F$, the result of Alexander and 
 Hirschowitz (\ref{1.1}) implies that $VSP(F,k)=\emptyset$, when 
 $k<{\frac 16}(n+4)(n+3)$.  In \ref{apolarity} we defined a map 
 $\alpha_{C}$ that associates an apolar cubic $n$-fold to an empty 
 codimension two linear section of a canonical curve $C$ of 
 genus $g=n+4$. The following theorem shows that cubic 
 $n$-folds in the image of this map
 are special with 
 respect to the 
 possible powersum presentations as soon as $n\geq 7$. 
 
 \begin{theorem}{\label{3.1}} If $F$ is a cubic $n$-fold apolar to 
 a general codimension two 
 linear section of a general canonical curve of genus $g=n+4$, then $VSP(F, 2n+4)\not=\emptyset$.  
 \end{theorem}
 \proof  This is immediate from \ref{2.4} since a canonical curve has nondegenerate 
 dual variety and satisfies $({\ast\ast})$.
\qed
 \begin{re}{\label{3.2}} By Hurwitz' formula, the degree of the 
 dual variety of a canonical curve is $6g-6$, so $VSP(F, 2n+4)$ 
 contains at least $6n+18$ points.  We do not know whether there are 
 more.
 \end{re}
 
  For $n\leq 3$, the general cubic is apolar to a section 
 of a canonical curve.  This fact can be used to describe 
 completely the powersum presentations of the cubic form (cf. \cite{RS}).

 For $n=3,4,5$ the general canonical curve of genus $g=n+4$ is a linear section of a 
 homogeneous space of dimension at least $6$ (cf. \cite{Muk}).  For 
 $n=4,5$, these homogeneous spaces of dimension $8$ (resp. $6$) have nondegenerate dual 
 varieties. For a cubic $4$-fold $F$ apolar to a general canonical 
 curve of genus $8$ proposition \ref{2.4} gives us a $4$-dimensional 
 component of $VSP(F, 10)$. It is shown in \cite{IR} that this is 
 in fact all 
 of $VSP(F, 10)$.  
 
  A general canonical curve of genus $9$ is a linear section of the 
  symplectic grassmannian $Sp(3)/U(3)\subset \Gr (3,6)$.
  For a cubic $5$-fold $F$ apolar to a canonical 
 curve of genus $9$, which is contained in a codimension two linear 
 section doubly tangent to $Sp(3)/U(3)$,  proposition 
 \ref{2.4} gives us a $2$-dimensional 
 subvariety of $VSP(F, 12)$.  On the other hand, for a general cubic $5$-fold $F$ 
 it follows from \ref{1.1} that $VSP(F, 12)$ is finite. 
 
 The general canonical curve of genus 
    $10$ is not a section of a $K3$-surface (cf. \cite{Muk}), so only 
    $\ref{3.1}$ 
    applies i.e. $VSP(F, 16)\not=\emptyset$, while already $VSP(F, 15)\not=\emptyset$ for a general cubic 
    $6$-fold $F$.

 \vspace{.2 in}
 
 Authors' addresses:
 \vspace{.2 in}
 
 Atanas Iliev\par 
 Institute of Mathematics, 
 Bulgarian Academy of Sciences,\par 
 Acad. G. Bonchev Str., 8,\par 
 1113 Sofia, Bulgaria.\par 

 e-mail: ailiev@math.bas.bg
 \vspace{.2 in}
 
 Kristian Ranestad\par 
 Matematisk Institutt, UiO,\par 
 P.B. 1053 Blindern,\par 
 N-0316 Oslo, Norway.\par 
 
 e-mail: ranestad@math.uio.no
 \end{document}